# Dynamical Systems Theory compared to Game Theory: The case of the Salamis's battle

**Konstantina Founta[1], Loukas Zachilas[1]**

[1]Dept. of Economics, University of Thessaly, Volos, Greece

Email: kfounta@uth.gr, zachilas@uth.gr

## 1. Abstract

In this paper, we present an innovative non–linear, discrete, dynamical system trying to model the historic battle of Salamis between Greeks and Persians. September 2020 marks the anniversary of the 2500 years that have passed since this famous naval battle which took place in late September 480 B.C. The suggested model describes very well the most effective strategic behavior between two participants during a battle (or in a war). Moreover, we compare the results of the Dynamical Systems analysis to Game Theory, considering this conflict as a "war game".

## 2. Introduction

In recent years, many researchers have studied the players' behavior either through Game theory or through Dynamical Systems. Some of the notable works are Archan and Sagar [2] who present a possible evolutionary game-theoretic interpretation of non-convergent outcomes. They highlight that the evolutionary game dynamics is not about optimizing (mathematically) the fitness of phenotypes, but it is the heterogeneity weighted fitness that must be considered. They mention that heterogeneity can be a measure of diversity in the population. In our research, this is described by the asymmetry in the conflict. In addition, Toupo, Strogatz, Cohen and Rand [3] present how important the role of the environment of the game is for the decision-makers. They suggest simulations of agents who make decisions using either automatic or controlled cognitive processing and who not only compete, as well as affect the environment of the game. Moreover, they propose a framework that could be applied in several domains beyond intertemporal choices, such as risky choice or cooperation in social dilemmas.

In other words, it's not all about the tactical behavior (aggressively or defensively), but also the impact of the location that the battle is taking place.

Furthermore, the well-known evolutionary game "Hawk – Dove" has been used in several scientific fields to describe the effects of different behavioral changes in populations. Some interesting applications are presented. Altman and Sagar [4] apply this game in a flock of birds modeling their behavior. In addition, Souza de Cursi [5] examines the applicability of Uncertainty Quantification (UQ) in this game under the uncertainty of both the reward and the cost of an injury to determine the mean evolution of the system. Lastly, Benndorf, Martínez-Martínez and Normann [6] investigate the equilibrium selection and they predict a dynamical bifurcation from symmetric mixed Nash equilibrium to asymmetric pure equilibria in the hawk – dove game, which depends on the frequency of interactions of the population.

Regarding the previous studies, it has been observed that there are no comparison results between dynamical analysis and game theory. The motivation of the present research refers to identify not only the connection of terms of these two scientific fields but also to apply this attempt in a battle.

The created model approaches short–term battles between two participants (players), where one is weaker than the other opponent. Also, the parameters (that we use in Eq. 1, see below) are the most crucial factors to highlight the optimal way to achieve a decisive victory. Below, the game hawk – dove and its results are presented.

One of the most representative games of Evolutionary Game Theory is the so-called game "Hawk – Dove", which was originally developed by Smith and Price [7] to describe animal conflicts and is quite similar to our attempt. There are two animals (or two players) fighting for the same resource. Each of them can behave either as a hawk (i.e., fight for the resource) or as a dove (i.e., abandon the resource before the conflict escalates into a fight). Individuals have a benefit B if they win and a cost C if lose.

If a Hawk meets a Hawk, they will fight and one of them will win the resource; the average payoff is *(B-C)/2*. If a Hawk meets a Dove, the Dove immediately withdraws, so the payoff of the Dove is *zero*, while the payoff of the Hawk is *B*. If a Dove meets a Dove, the one who first gets hold of the resource keeps it, while the other does not fight for it; average payoff *B/2*. The strategic form of the game is given by the payoff matrix (1):

$$Payoff_{H,D} = \begin{pmatrix} (B-C)/2 & B \\ 0 & B/2 \end{pmatrix} \quad (1)$$

## 3. Results of the game "Hawk – Dove"

We set random values in the benefit $B = 2$ if a player wins, and in the cost $C = 1$ if a player loses in the payoff matrix (1). Using the Gambit[a] software (16.0.1), we find Nash equilibriums and the dominant strategy.

| | HAWK | | DOVE | |
|---|---|---|---|---|
| **Player 1** Payoff: 2 | | | | |
| HAWK | ½ | 0 | 2 | 0 |
| **Player 2** Payoff: 0 | | | | |
| DOVE | 0 | 2 | 1 | 1 |

Profiles 1 — All equilibria by enumeration of mixed strategies in strategic game

| # | 1: HAWK | 1: DOVE | 2: HAWK | 2: DOVE |
|---|---|---|---|---|
| 1 | 1 | 0 | 1 | 0 |
| 2 | 1 | 0 | 0 | 1 |

*Figure 1. The results of "Hawk - Dove" game.*

Figure 1 shows us the payoff matrix and the two Nash equilibriums. If both players behave as a Hawk, the one who first injures the other wins. We set the player 1 starts and injures the player 2, thus player 1 wins. If someone behaves as a Hawk and the other behaves as a Dove, then the player with the aggressive behavior (Hawk) wins and takes all the resource. If both players behave as a Dove, then they share the resource.

Regarding Nash equilibriums, there are two pure strategies. On the one hand, both players behave as Hawks and on the other hand, player 1 behaves as a Hawk and player 2 as a Dove. Moreover, we can observe that player 1 behaves as a Hawk in both cases and player 2 behaves either as a Hawk or as a Dove, but in each case player 1 wins.

We should note that if player 2 injures first player 1, the Nash equilibriums would be different.

---

[a] McKelvey, Richard D., McLennan, Andrew M., & Turocy, T. L., (2014).

*Figure 2. Dominant Strategy.*

Figure 2 shows the dominant strategy of the game, where player 1 behaves as a Hawk independently of the player's 2 behavior (i.e., Hawk or Dove). Therefore, the first dominant strategy may not be effective, because both players behave as Hawks and player 1 wins the half of the resource and not maximize his profit. Although, if the player behaves as a Hawk, knowing that the other player behaves as a Dove, then he takes all the resource (maximum profit). Thus, we believe that the second Nash equilibrium is more effective and optimum strategy.

## 4. The Dynamical Model

It is widely acknowledged that the military strategy is the combination of «*ends, ways and means*» [8]. In our attempt to study the strategic behavior of two warring parties, we developed an innovative non-linear discrete system of two equations based on the above phrase. The main objective of the model is to simulate the way by which the two opponents behave strategically, where the one is weaker than the other.

Simultaneously, in Game Theory, the war is considered as a dynamic game where the strategies of the players are studied by calculating their optimal strategy (Nash equilibrium). In the present paper, we compare the results of the Game Theory to those from the analysis of the discrete dynamical system. At the end of the analysis, the optimum and effective strategy for both participants (players) will be suggested.

The model, which is applied in short-term conflicts and describes the strategic behavior of each participant, is given by Eq. 2:

$$\begin{cases} x_{t+1} = P_x + TN_x - G \cdot (D_y + E_x) \cdot 4 \cdot y_t \cdot (1 - y_t) \\ y_{t+1} = P_y + TN_y - (1 - G) \cdot (D_x + E_y) \cdot 4 x_t \cdot (1 - x_t) \end{cases} \quad (2)$$

where:

$x_t$: The strategic behavior of the participant (player) *x* at the time *t*.

$y_t$: The strategic behavior of the participant (player) *y* at the time *t*.

$x_{t+1}$: The optimal strategic behavior of the participant (player) *x* at the (next moment of) time *t* + 1.

$y_{t+1}$: The optimal strategic behavior of the participant (player) *y* at the (next moment of) time *t* + 1.

We consider $x_t, y_t, x_{t+1}, y_{t+1} \in [0,1]$, because the logistic equation is defined in [0,1], which is derived from the study of biological populations reproduced in discrete time [9]. It's the evolution of the population model of Malthus [10] and shows that the exponential growth cannot tend to infinity, but there is a critical point, i.e., a saturation. In other words, it is not possible for someone to win and the other to continuously lose. Also, each optimal strategic behavior, at the time *t*, affects the next move – strategic behavior, at the time *t* + 1, of the opponent.

In addition, we can interpret the values of variables (and parameters, as shown below) as percentages or probabilities, which help us to explain the results; these are also explained through the Game Theory.

Moreover, if the value of $x_{t+1}$ (or $y_{t+1}$, respectively) equals to 0, it indicates the fully defensive strategic behavior of participant *x* (or *y* respectively), while if it equals to 1, then it indicates the fully aggressive behavior of participant *x* (or *y* respectively).

The parameters of Eq. 1 are the main and most important factors that could affect the strategic behavior of *x* (or *y*, respectively). In particular:

The parameter $\boldsymbol{P_x}$ represents the strength (economic, military, population, territorial) of *x* and $\boldsymbol{P_y}$ is the strength of *y*, respectively. These two parameters indicate the substance of each form of social organization compared to the other.

$TN_x$ and $TN_y$ represents the Technological Naval capability and evolution of *x* and *y* respectively. These two parameters are also defined in comparison with the technological capability and evolution of the other participant and describe the *means* mentioned by Lykke [8].

The parameter *G* represents the geographical location (geophysical terrain) of the area where the battle or the war is taking place. We believe that this is another part of the military strategy, namely the *ways* [8]. Trying to emphasize the importance of this parameter and how it can be an advantage or disadvantage for each participant, we set in the first equation as *G* and in the second equation as **1 − *G***. The closer to the 1 the value of the parameter, the easier the geophysical terrain of the area is.

The parameter $D_x$ represents the damages caused by *x* to *y* and respectively, $D_y$ represents the damages that *y* brings to *x*. The damages which we refer to may be economic, territorial, military, etc. or even deception and damaging of the psychological part of the opponent. Moreover, these two parameters complete the last part of the military strategy, namely the *ends* [8].

The parameter $E_x$ represents the expenses of participant *x* and $E_y$ the expenses of participant *y*, respectively. In other words, these denote the preparation costs of each participant for a battle (or war), compared to each other.

All the parameters that have been presented above should belong to [0,1]. Namely, $P_x, P_y, TN_x, TN_y, G, D_x, D_y, E_x, E_y \in [0,1]$.

In the next section, we present the dynamic analysis and the results from the application of Eq. 1 in naval battle of Salamis.

## 5. The case of (naval) Battle of Salamis

The naval battle of Salamis was an important battle of the second Persian invasion in Greece and has been estimated to being held on September 28[th], 480 BC in the Salamis straits (in the Saronic Gulf near Athens). The two warring parties were the Greeks (Hellenic alliance) and the Persian Empire [11].

After the fall of Thermopylae, the Persians went ahead to Athens. The Greeks had been advised by the Oracle of Delphi, that only the "*wooden walls*" would save them, and they considered that this referred to a fight in the sea [12].

A few days before the battle, the war council of the Greek admirals had to decide the geographic location of the battle. On the one hand, the Spartan General Evriviades proposed to fight at the Isthmus of Corinth, under the main argument that in case of failure it would be possible for them to continue to fight into the center of the Peloponnese. On the other hand, the Athenian General Themistocles insisted to fight in Salamis's straits. He believed that if he forced the Persians to attack there, the numerous Persian ships couldn't extent highlighting their dominance. Ultimately, the council considered that Themistocles' argument was better and decided to support it [13].

The Greek fleet was estimated by Herodotus in 380 triremes and Aeschylus gave a round 300 triremes, but we can't be certain for the exact number. On contrary, the Persian fleet was estimated in 500-600 triremes[b]. Herodotus describes the Persian ships as "*better sailing*", when compared to the Greek fleet. This may be attributable to a combination of factors such as lightness of materials and structure of the ship, better seamanship, and more extensive naval experience. The triremes of Hellenic alliance were heavier and more durable. However, Herodotus reports that these ships were equipped with an embolism, with which they sank the enemy ships. They used two attacking maneuvers: *diekplous*, (i.e., attack from the rear or sides with a sharp turn) and *periplous*, (flanking or enveloping move, which generally gave an extra benefit against superior numbers in open water). The purpose of both was to ram the enemy in the side. In this way, they achieved serious damages or even the complete destruction of the Persians ships. On the contrary, the Persian tactic was "ramming and boarding" [14].

---

[b] Aeschylus, writing decades earlier, also gives 1,207 triremes, but Herodotus writes, shortly before battle took place, that the Persian fleet wasn't much bigger than Greek. Because of a weather phenomenon (storms) 600 ships sank (400 at the coast of Magnesia, north of Artemisium and 200 in Euboea).

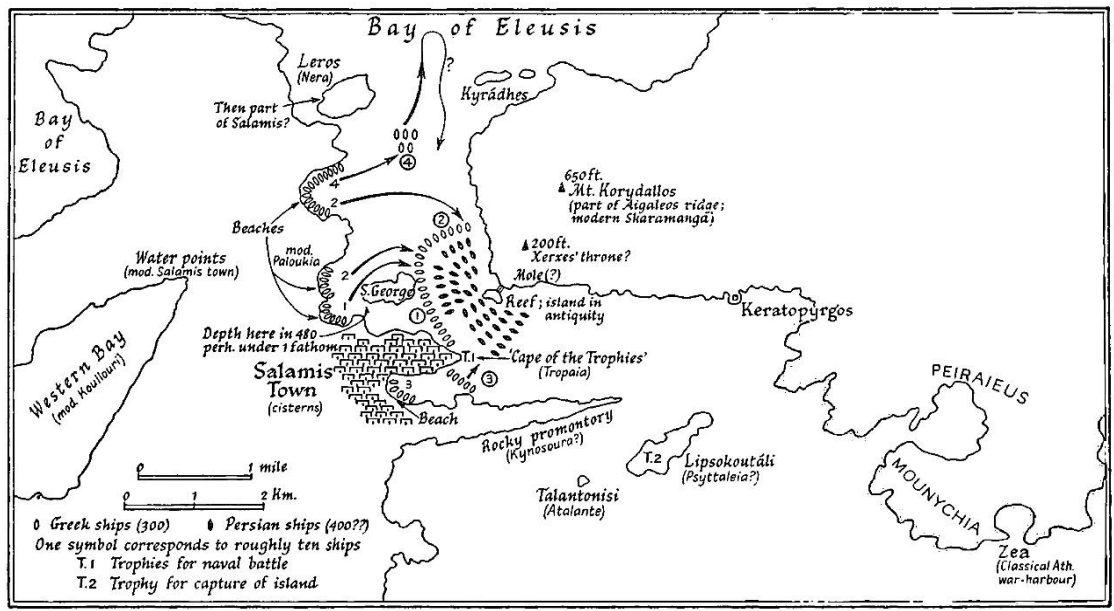

*Figure 3. The battle of Salamis. Source: Burn, A. R. (1962). Persia and the Greeks, New York: Minerva Press, p.452*

At the dawn (if the date of the battle was indeed 28[th] Sept.), the two fleets were ready for the naval conflict. Xerxes, sure of his victory, placed his throne on mountain Aigaleo (see Fig. 3), to enjoy the war spectacle. The narrowness of the space and the limited extent of the sea did not allow the Persians to use the major of their force in the front line. Thus, the number of ships was approximately equal. In this naval battle, the bravery and dexterity of the Greek fleet played an important role. They fought aggressively to defend their moral values and their freedom [13].

Herodotus reports that "*the Greeks fought with discipline and held their formation, but the Persians did not seem to be following any plan, so things were bound to turn out for them as they did*". Also, Aeschylus mentions that Themistocles must be given the credit for their battle and the winning tactics. The turning point of the battle came as the Persians "*suffered their greatest losses when the ships in their front line were put to fight and those following, pressing forward to impress the King* (i.e., Xerxes) *with their deeds, became entangled with them as they tried to escape*", as Herodotus comments [13].

The naval battle evolved rapidly and by the noon it was visible that the Greeks would win. The Persian fleet had crushed, while the Greek fleet continued to haunt them, killing the helpless, non–swimming soldiers. This brought the battle to an end, leaving the Greek force in full control of the straits [14].

When the battle was over, a Roman source mentions that Greeks lost more than 40 triremes and Persians more than 200 ones [14]. The victory of the Greek force was of major importance, since they managed to cause the collapse of the Persian morale, which is evidenced by the abandonment of the battle. In addition, the right decision of Themistocles for the geographic location of the naval battle was one of the most intelligent movements to bring the Greek victory.

## 6. Applying the model in naval battle of Salamis – Approaching the reality

Starting the dynamical analysis of the naval battle of Salamis, we set the initial conditions (in Eq. 1), which represent as much as possible the historical events of the battle. Specifically:

(a) We set Greeks as the weak participant – player (x) and Persians as the powerful participant – player (y).

(b) The strength of Hellenic alliance, $P_x = 0.25$ and the strength of Persian empire, $P_y = 0.8$. The values of parameters describe the triremes (in quality and quantity) of each fleet. As mentioned above, the Greek fleet was estimated 300 – 380 triremes and the Persian fleet was estimated in 500 – 600 triremes.

(c) The technological naval capability of Greeks, $TN_x = 0.7$ and the technological naval capability of Persians, $TN_y = 0.35$. According to historical documents, the Greek fleet were well-trained in relation with Persian sailors. Moreover, the Greek triremes had an embolism in front to attack the enemies' ships.

(d) The geographic location of the naval battle, $G = 0.4$, i.e., the Salamis straits, which are an advantage point for the Greek fleet. In a very recent paper, Zerefos et al. [15] study the climatically prevailing weather conditions during the battle of Salamis. They mention that Themistocles was aware of the wind patterns on the day of the battle and the knowledge of the local wind climatology in combination with the narrow of the location must have been a critical argument to confront the Persian fleet in the straits of Salamis. The variation of the wind in the Saronic Gulf was west – northwest at previous night and early in the morning being replaced by a southeast wind after 10.00 A.M. and through early evening. This information was very useful to the strategic plan of the Greek fleet to trap the Persians in Salamis and led to one of the greatest victories in history.

(e) The damage caused to Persian side was huge, so we set $D_x = 0.8$ and $D_y = 0.2$. As mentioned above, Greeks lost more than 40 triremes and Persians more than 200 ones.

(f) The preparation costs of this battle for each participant: $E_x = 0.3$, $E_y = 0.7$, respectively. According to Kyriazis and Zouboulakis [16], 100 new Athenian triremes were built under the Athenian Naval Law of Themistocles. Each one cost one talent (6000 ancient drachmae), so the total cost was 100 talents (or 600.000 ancient drachmae). In 480 BC, the Athenian fleet was comprised of 200 triremes, equivalent to the two thirds of the total Greek strength. However, the Persian ships were similar in shape, so we assume that the cost of each ship was similar. Thus, it is obvious that the Persians spent more money to support their expedition to the Greek territories than the Greeks.

With these initial conditions, we solve the system (Eq.1), by using the mathematical software Maxima[c] (5.39.0), calculating the equilibrium points. Then, we study more extensively the behavior of the model, and we present bifurcation diagrams and timeseries diagrams using the software E&F Chaos[d].

Solving the system (Eq.1) there are two equilibrium points: $E_1$ ($x^* = 0.75$, $y^* = 0.475$) and $E_2$ ($x^{**} = 0.96$, $y^{**} = 1.012$). According to Game Theory, these two fixed points are considered as Nash Equilibriums [1]. Below, the stability of the fixed points will be examined.

The Jacobian matrix is:

$$J = \begin{pmatrix} 0 & 0.8y - 0.8(1-y) \\ 3.6x - 3.6(1-x) & 0 \end{pmatrix} \quad (3)$$

We calculate the Jacobian matrix at the equilibrium point $E_1$:

$$J^* = \begin{pmatrix} 0 & 0.038 \\ 1.803 & 0 \end{pmatrix} \quad (4)$$

The determinant of $J^*$ is det ($J^*$) = 0.069 > 0.

---

[c] https://sourceforge.net/projects/maxima/files/Maxima-Windows/5.39.0-Windows/
[d] E & F Chaos: written by Diks, C., Hommes, C., Panchenko, V., van der Weide, R., (2008).

The trace of $J^*$ is trace($J^*$) = 0.

The eigenvalues of $J^*$ is $(0.264i, -0.264i)$; two complex roots.

The discriminant is $\Delta = trace(J^*)^2 - 4 \cdot det(J^*) = -0.2788 < 0$.

Therefore, the equilibrium point $E_1$ is a <u>stable – center</u>.

Studying the second fixed point $E_2$, the Jacobian matrix at the equilibrium point is:

$$J^{**} = \begin{pmatrix} 0 & 0.082 \\ 3.314 & 0 \end{pmatrix} \quad (5)$$

The determinant of $J^{**}$ is det ($J^{**}$) = $-2.718 < 0$.

The trace of $J^{**}$ is trace($J^{**}$) = 0.

The eigenvalues of $J^{**}$ is $(1.648, -1.6487)$; two real roots.

The discriminant is $\Delta = trace(J^{**})^2 - 4 \cdot det(J^{**}) = 10.874 > 0$.

Therefore, the equilibrium point $E_2$ is a <u>saddle point</u>.

Consequently, we accept the fixed point $E_1$ ($x^* = 0.75$, $y^* = 0.475$) and reject $E_2$ ($x^{**} = 0.96$, $y^{**} = 1.012$), because the value of $y^{**}$ is greater than 1.

Thus, we continue the analysis for the fixed point $E_1$. Interpreting this equilibrium point, we confirm the aggressive (strategic) behavior of Greeks; since the value of $x^*$ is close to 1 and the mild (strategic) behavior of Persians; since they thought it would be an "*easy win*".

Indeed (historically), the courage of the Greeks, their technological naval skills, and the advantageous geographical location contributed to this aggressive behavior. As far as the Persians are concerned, their mild (strategic) behavior is due to the fact that they underestimated their enemy, since they regarded that the Greeks are an easy target, and they would achieve a decisive victory.

Connecting the game "Hawk – Dove" to the naval battle of Salamis, player 1 (red) is "Persians" and player 2 (blue) is "Greeks" (Fig. 4). The Hellenic alliance had an aggressive behavior (Hawk) and the Persians behaved as a Dove. According to the Nash equilibriums that have been mentioned above (See 3), the Greeks (player 2) should behave as a Hawk (i.e., aggressive), regardless of Persian's behavior, to win this battle.

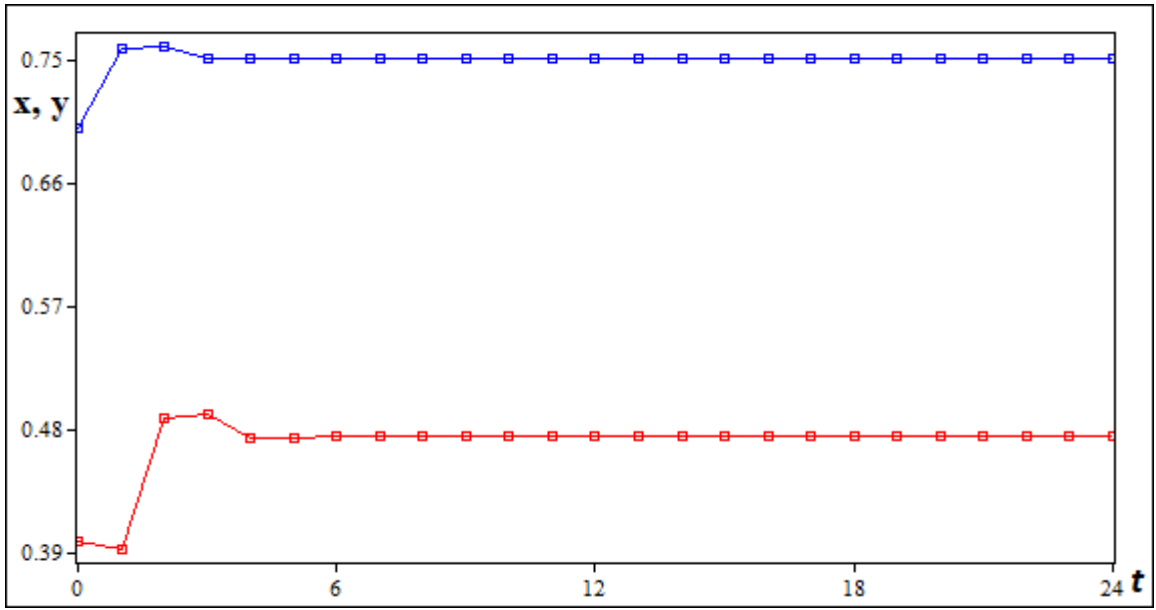

*Figure 4. Time series diagram - x (Greeks; blue) and y (Persians; red) – t(hours).*

Figure 4 shows us how the two warring parties behave (strategically). It represents the optimal strategic behavior of Greeks and Persians in Salamis straits for a time interval of 24 hours. We can observe an oscillation, at the beginning, until $t = 6$ h. (both lines) and then it is normalized and balanced. That means that the duration of the main battle was approximately 6 hours. Indeed, according to historical documents, the battle started at dawn (approximately at 06:00 am) and the Greek victory was visible at noon.

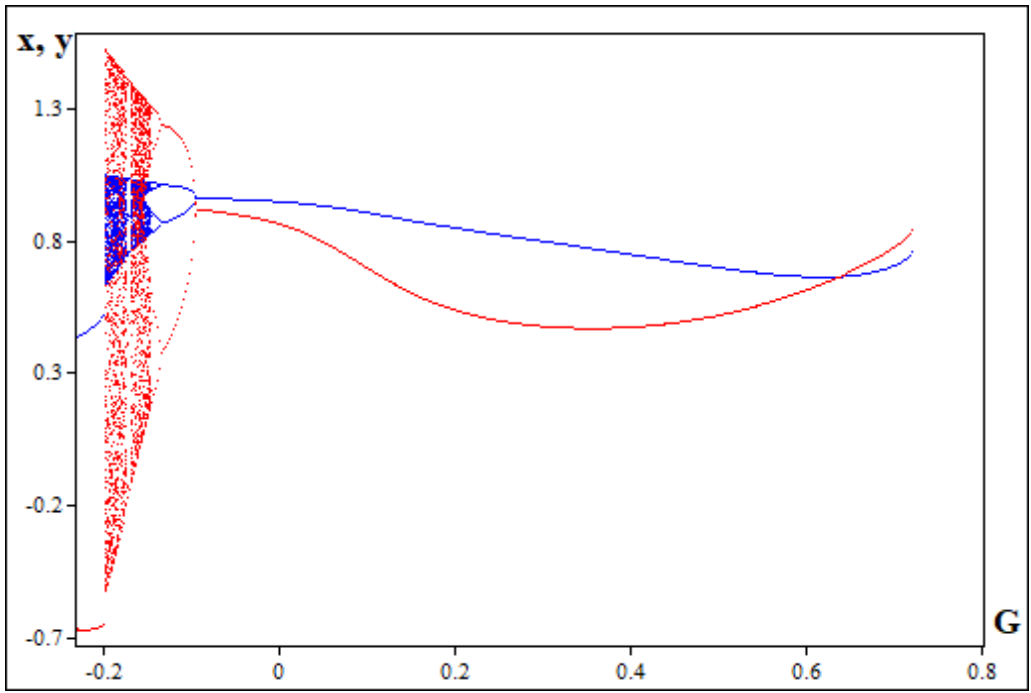

*Figure 5. Bifurcation Diagram for different values of G – x (Greeks; blue) and y (Persians; red).*

*x, y – vertical axis and G – horizontal axis.*

Figure 5 presents the strategic behavior of Greeks (blue) and Persians (red) as the parameter $G$ changes. We observe for the positive values of $G$, the blue line is above the red until $G = 0.64$ (critical value) and for $G > 0.64$ the red line is above the blue. The increase of the value of parameter signifies the change of the geographical location (a more open sea), which becomes more difficult for Greeks and in contrary easier for Persian. Thus, we approve that if the location of the naval battle was in an open sea, the Persians would have a crucial advantage, which would possibly lead to win this conflict.

Although we did not study the negative values of parameter $G$, we believe that there are some unpredictable geophysical factors (e.g., meteorological phenomena to influence the outcome of the conflict), which are surprisingly interesting. Specifically, we refer to weather conditions, such as air, ripple, etc., which can affect the geophysical terrain of the area. Due to these weather phenomena, period doubling bifurcations and chaos appear and we cannot predict what could happen in the battle for these values of $G$.

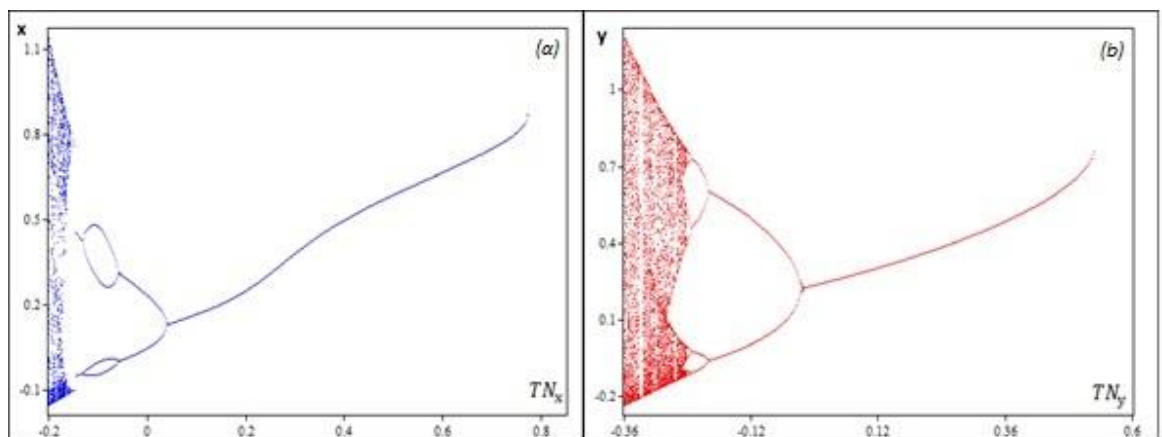

*Figure 6. Bifurcation diagrams for different values of parameters $TN_x$ and $TN_y$.*

*Fig. 6a: x (Greeks; blue) – vertical axis and $TN_x$ – horizontal axis,*

*Fig. 6b: y (Persians; red) – vertical axis and $TN_y$ – horizontal axis.*

Figure 6 depicts the technological evolution and capability of $x$ (Greeks; blue) and $y$ (Persians; red), respectively. In the left diagram (Fig. 6a), for the negative values of parameter $TN_x$, we can distinguish a pair of bubble bifurcations, while afterwards we have the well-known period-doubling scenario to chaos. A possible interpretation of this chaotic scenario is the uncertainty of Greeks in technological capability – first attempts to construct ships. The first ships, as Krasanakis [17] mentions, were floating planks and carved tree trucks only with oars. Since the ships were primitive, the situation was unstable (there is chaos in this range of values) because they were not capable to fight in naval battles. Later, the sails were invented, which gave high speed to ships, and they were consisted no more than wood but iron. For this reason, we have bubble bifurcations, which indicate the technological alternatives that existed for the construction of the ships. In the interval of positive values of $TN_x$, there is stability with two fixed points. Here, it is the beginning of better shipbuilding ability and new expertise ship construction. Finally, there is one equilibrium point which shows the better version of ships, of that period, namely Triremes. Triremes were wooden warships which move either with sails or oars. Moreover, in the positive values of the parameter $TN_x$, the increasing of the slope of the curve is visible, which, on the one hand, it means that in 480 BC the triremes were an innovation in shipbuilding and on the other hand, it shows the excellent naval capability of the Greeks.

Persians, through the years, developed technological equipment because of their expansive mania to conquer Greece. Comparing the Figures 6a and 6b, it seems that Persians had a lower technological development than Greeks, since they focused more on land army than on warships. Their ships were mainly used as troopships rather than battleships [14].

## 7. "What if…": Two scenarios by changing the geographical conditions of the battle

In this section, we will study two alternative scenarios, that were indeed discussed in Greek Generals' meeting. In the first scenario we changed the place of the naval battle, which now is supposed to take place in a more open sea, (i.e., openly the Saronic Gulf), while in the second scenario we changed the place to a mixed battle (part of the battle takes place in Saronic straits, and part at the mainland of Isthmus of Corinth).

**(a) Openly the Saronic Gulf**

In Eq. 1 we keep the values of the parameters constant and change only the value of parameter *G* to 0.64. The increase of this parameter changes the geophysical landscape to a more open sea. Now the Greeks are not in an advantaged position and the Persians have a more dominant position, because they have the chance to add more ships in the battle.

We solve Eq. 1, and we get two equilibrium points: $E_1$ ($x^* = 0.667$, $y^* = 0.67$) and $E_2$ ($x^{**} = 0.904$, $y^{**} = 0.963$). According to Game Theory, these fixed points are Nash equilibriums. The stability of these points is studied, again, by the Jacobian matrix:

$$J = \begin{pmatrix} 0 & 1.28y - 1.28(1-y) \\ 2.16x - 2.16(1-x) & 0 \end{pmatrix} \quad (6)$$

We calculate the Jacobian matrix at the fixed point $E_1$:

$$J^* = \begin{pmatrix} 0 & 0.436 \\ 0.722 & 0 \end{pmatrix} \quad (7)$$

The determinant of $J^*$ is det $(J^*) = -0.314 < 0$.

The trace of $J^*$ is trace $(J^*) = 0$.

The discriminant is $\Delta = trace(J^*)^2 - 4 \cdot det(J^*) = 1.259 > 0$.

Therefore, the equilibrium point $E_1$ is a <u>saddle point</u>.

Studying the second fixed point $E_2$, the Jacobian matrix at the equilibrium point is:

$$J^{**} = \begin{pmatrix} 0 & 1.184 \\ 1.745 & 0 \end{pmatrix} \quad (8)$$

The determinant of $J^{**}$ is det $(J^{**}) = -2.068 < 0$.

The trace of $J^{**}$ is trace $(J^{**}) = 0$.

The discriminant is $\Delta = trace(J^{**})^2 - 4 \cdot det(J^{**}) = 8.273 > 0$.

Therefore, the equilibrium point $E_2$ is a <u>saddle point</u>.

For the first fixed point ($x^* = 0.667$, $y^* = 0.67$), we notice that the values of $x^*$ and $y^*$ are very close. This means that $G = 0.64$ is close to a critical value. The "new" geographical location (a more open sea, than Salamis's straits) gives to the Persians the advantage to include more warships in the naval battle. The nature of this fixed point (saddle point) leads our thought that this would be a turning point changing the whole outcome of the conflict, while the increase of the number of Persian ships makes the situation unstable.

For the second fixed point ($x^{**} = 0.904$, $y^{**} = 0.963$), we notice that the values of $x^{**}$ and $y^{**}$ are – again – very close, while very close to 1. This means that both opponents have a very aggressive behavior (we may assume that it would be a conflict between two Hawks, according to game "Hawk – Dove"). This is a non-effective scenario (strategy) for Greeks, because the Greek's benefit of such a conflict would be lower than the scenario: "Greeks / Hawk" versus "Persians / Dove". Nevertheless, it can be classified as an unstable situation for the same reasons to the first fixed point.

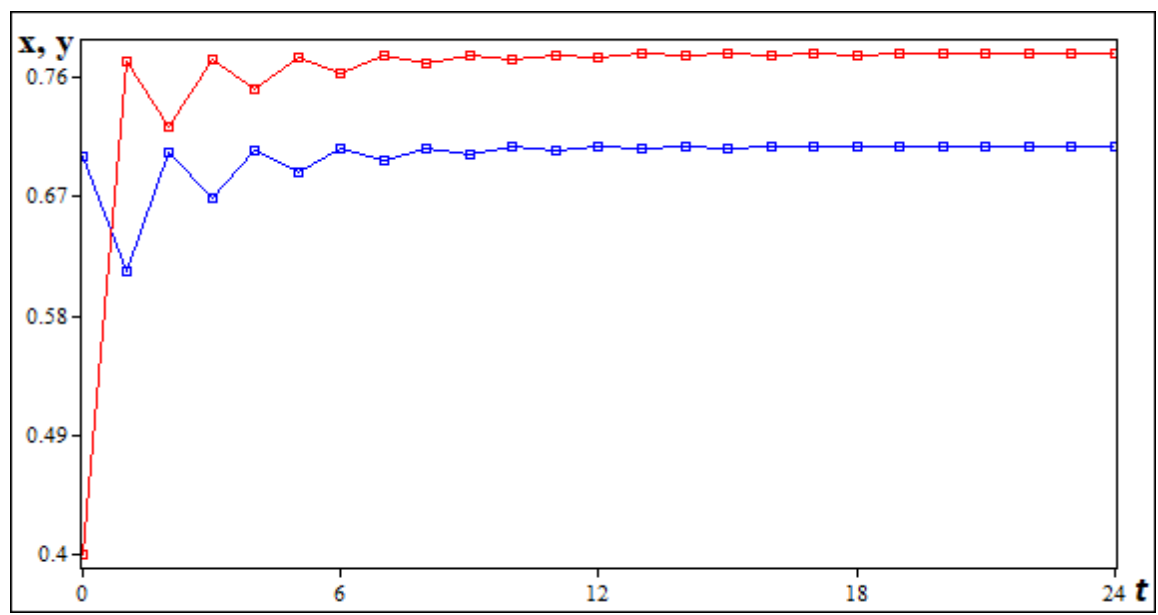

*Figure 7. Time series diagram - x (Greeks; blue) and y (Persians; red) in the open Saronic Gulf – t(hours).*

Figure 7 presents the strategic behavior of Greeks and Persians in the open Saronic Gulf. At first, we observe an oscillation up to $t = 10\sim11$ h (both lines), while later it is normalized and balanced. This means that the duration of the battle would be approximately 11 hours (twice as much than the real duration) to determine the outcome of the conflict. Moreover (and this is the real interesting outcome) the red line is above the blue one, which means that due to the open sea the Persians could increase the number of their ships and thus winning the battle.

**(b) Isthmus of Corinth**

In Eq. 1 we keep the values of the parameters constant, and we change only the value of parameter $G$. Now, we set $G = 0.7$. Increasing more the value of parameter $G$, we try to present the mixed battle. Part of the battle takes place in Saronic straits, and part at the mainland of Isthmus of Corinth, where Persians could add even more ships in the naval battle. This is the geographical location that was proposed by Spartan General Evriviades.

We solve Eq. 1, and we get (again) two equilibrium points: $E_1$ ($x^* = 0.7$, $y^* = 0.77$) and $E_2$ ($x^{**} = 0.85$, $y^{**} = 0.92$). According to Game Theory, these fixed points are Nash equilibriums. We have studied (again) the stability of these points by calculating the Jacobian matrix:

$$J = \begin{pmatrix} 0 & 1.4y - 1.4(1-y) \\ 1.8x - 1.8(1-x) & 0 \end{pmatrix} \quad (9)$$

The Jacobian matrix at the fixed point $E_1$:

$$J^* = \begin{pmatrix} 0 & 0.778 \\ 0.708 & 0 \end{pmatrix} \quad (10)$$

The determinant of $J^*$ is det $(J^*) = -0.583 < 0$.

The trace of $J^*$ is trace $(J^*) = 0$.

The discriminant is $\Delta = trace(J^*)^2 - 4 \cdot det(J^*) = 2.334 > 0$.

Therefore, the equilibrium point $E_1$ is a <u>saddle point</u>.

Studying the second fixed point $E_2$ ($x^{**}= 0.85$, $y^{**}= 0.92$), the Jacobian matrix at the equilibrium point is:

$$J^{**} = \begin{pmatrix} 0 & 1.194 \\ 1.277 & 0 \end{pmatrix} \quad (11)$$

The determinant of $J^{**}$ is det $(J^{**}) = -1.526 < 0$.

The trace of $J^{**}$ is trace $(J^{**}) = 0$.

The discriminant is $\Delta = trace(J^{**})^2 - 4 \cdot det(J^{**}) = 6.104 > 0$.

Therefore, the equilibrium point $E_2$ is also a <u>saddle point</u>.

We can observe for the first fixed point $E_1$ that the values of $x^*$ and $y^*$ are very close, while, for the second fixed point $E_2$ the values of $x^{**}$ and $y^{**}$ have a slight deviation. In both fixed points, the value of $y$ is greater than $x$, which means that in both cases the Persian fleet could win this conflict anyway. Therefore, since in both cases the fixed points are saddle ones, the situation is unstable and Persian dominance in this mixed battle is indisputable.

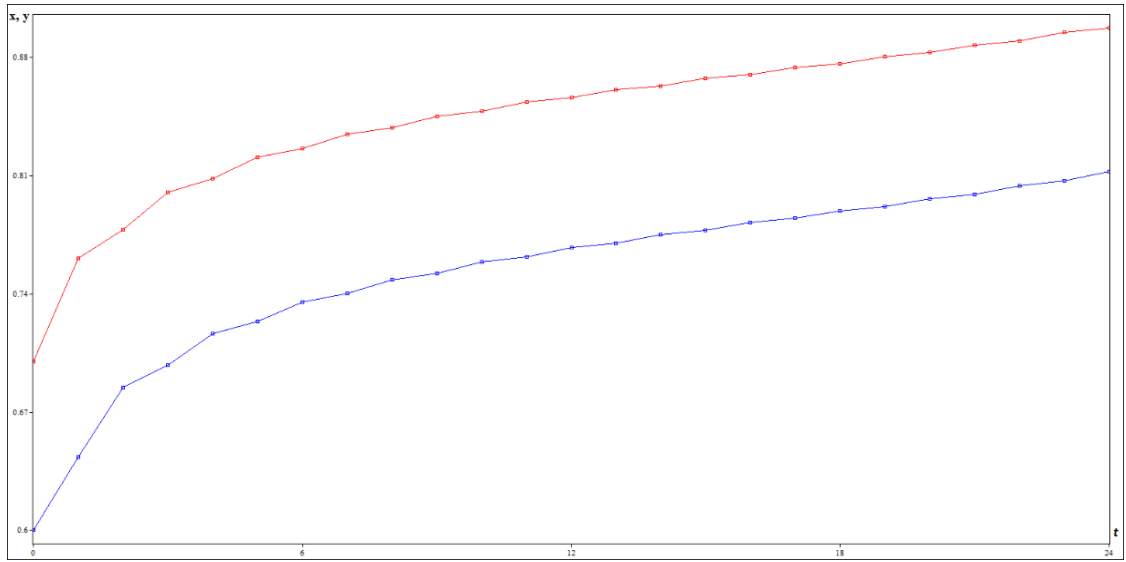

*Figure 8. Time series diagram - x (Greeks; blue) and y (Persians; red) in the second scenario – t(hours).*

Figure 8 represents the strategic behavior of Persians and Greeks in the mixed battle (Salamis's straits and Isthmus of Corinth) throughout time. We can see that both sides have a continuous upward trend, and the red line (Persians) is above the blue line (Greeks) for all the time interval, which means that the Persian's victory would have taken place from the outset. Moreover, it is obvious that this geographical location would be the worst choice for the Greeks and Evriviades's view would have led the Greeks to a crushing defeat!

## 8. Studying bilateral damages of two opponents

In this section, we present the results of three hypothetical scenarios concerning the damages caused by Greeks to Persians and vice versa. We will keep all the rest parameters at their initial values, representing that the battle took place in Salamis's straits.

In the first scenario, we assume that $D_x = 0.5$ and $D_y = 0.5$, i.e., the damage that Greeks caused to Persians is 50% (of their total armament) and vice versa. In the second scenario, we set $D_x = 0.3$ and $D_y = 0.7$, i.e., Greeks cause 30% damage to Persians and Persians cause 70% damage to Greeks, respectively. In the last scenario, we assume that $D_x = 0.8$ and $D_y = 0.2$, i.e., Greeks cause 80% damage, while the Persians cause 20% damage to Greeks. The time series diagrams for each scenario are presented in Figures 9, 10 and 11.

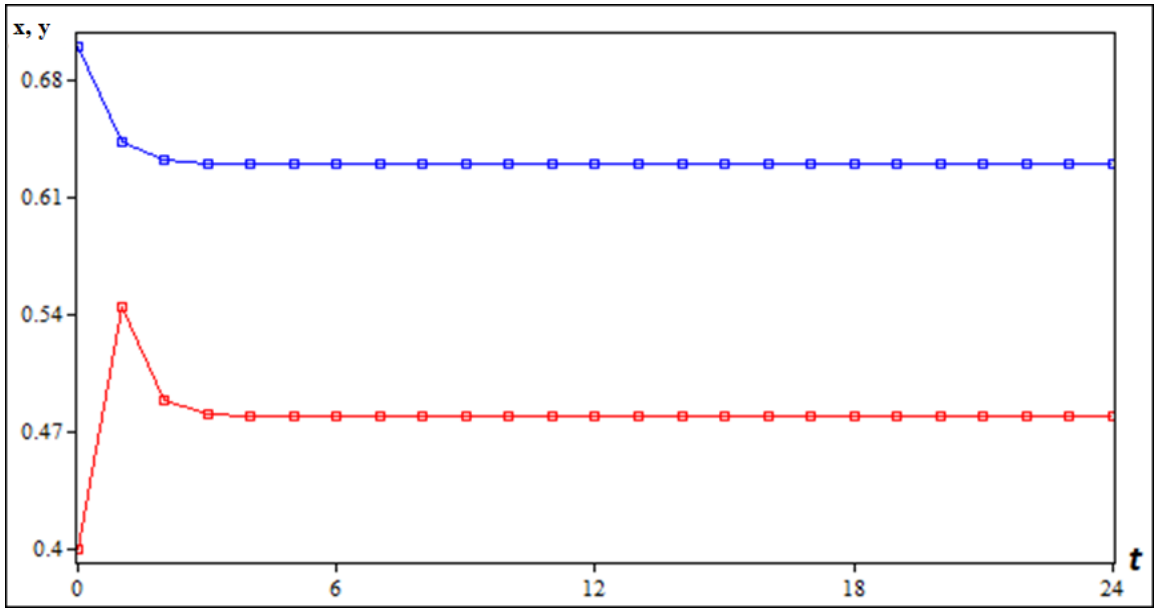

*Figure 9. Scenario 1: Time series diagram - x (Greeks; blue) and y (Persians; red) – t(hours).*

Figure 9 shows the effect of bilateral damages on the strategic behavior of Greeks (blue) and Persians (red), respectively. The supremacy of the Greeks is evident from the beginning, since, on the one hand, they had a better technological ability and an advantageous geographical location and, on the other hand, the fact that 50% of the damage to the opponent would be capable of bringing the Greeks a decisive victory.

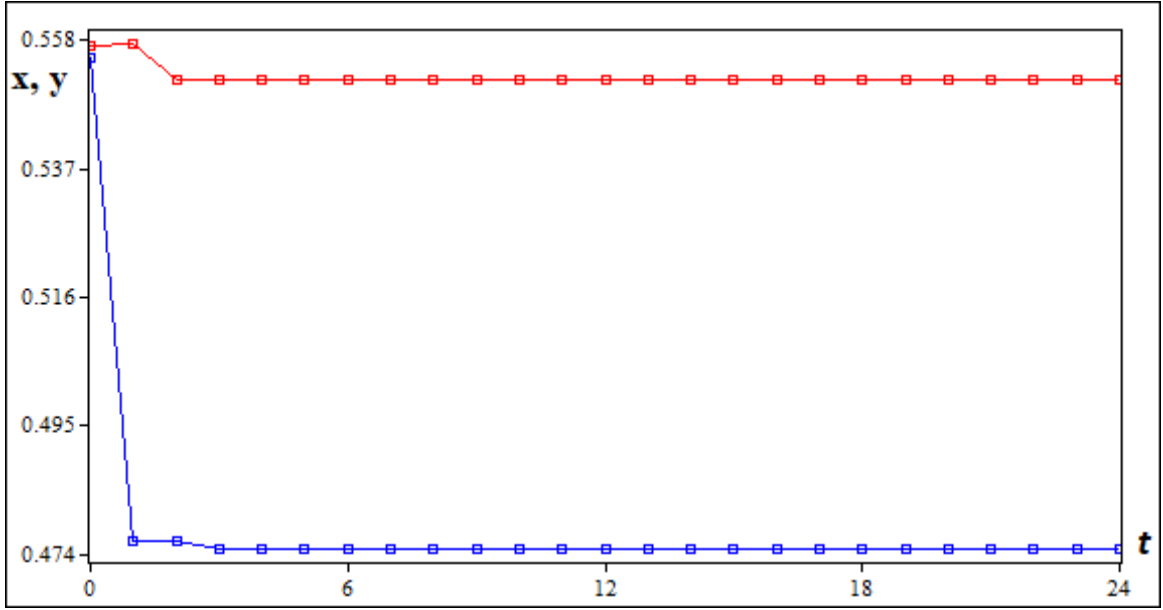

*Figure 10. Scenario 2: Time series diagram - x (Greeks; blue) and y (Persians; red) – t(hours).*

Figure 10 shows the effect of bilateral damages on the strategic behavior of Greeks (blue) and Persians (red), respectively. If Persians had caused more damage to the Greeks, it would be a reason to win this naval battle in a very short time. If Persians had caused 70% damages to Greeks, we are talking about a total destruction of Hellenic alliance and a remarkable victory by the Persians.

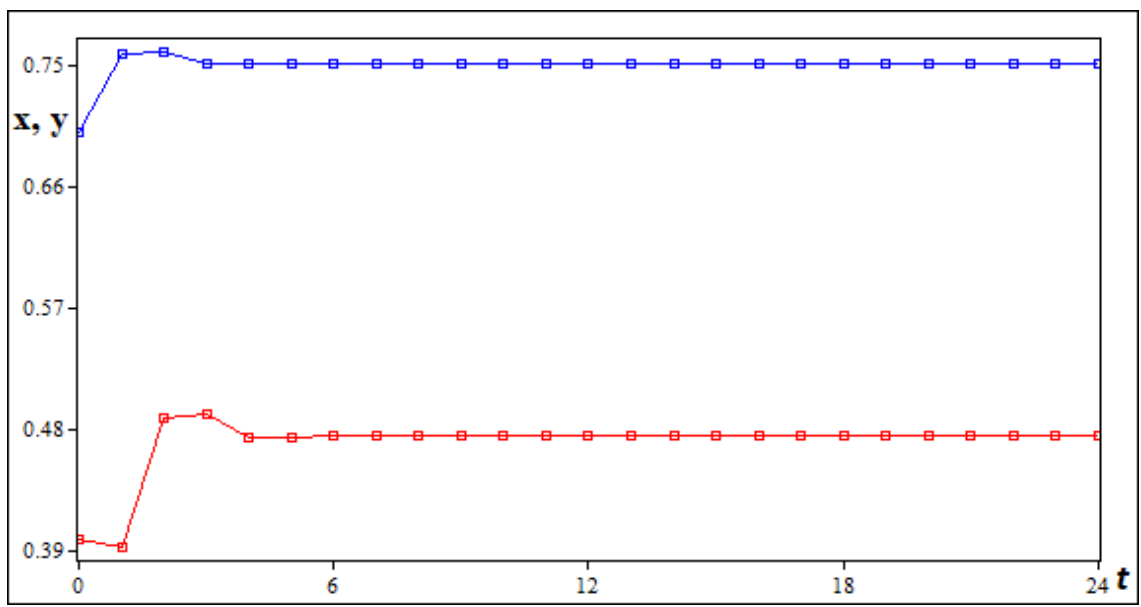

*Figure 11. Scenario 3: Time series diagram - x (Greeks; blue) and y (Persians; red) – t(hours).*

Figure 11 shows the real historical evolution of events, since the values of parameters realistically approach what happened in his naval battle. The supremacy of the Greeks, at the beginning, in the Salamis straits is owing to naval tactic "*diekplous*" and the advantageous geographical location. In this way, Greeks achieved the decisive victory against to Persian empire.

## 9. Conclusions

In this paper, an innovative non – linear discrete model has been presented, which simulates the optimum strategic behavior of two warring parties for short–term battles. In addition, we try to compare this model with the classical Game Theory, applying this attempt in the naval battle of Salamis. Based on the results we have extracted, we (mathematically) proved the historical events of this conflict, and we concurrently studied some alternative hypothetical scenarios. Specifically, by changing the geographical location of the conflict, we prove that the optimum location for Greeks was the Salamis straits and, on the contrary, the worst location for achieving the decisive victory would be Isthmus of Corinth. Moreover, we study various scenarios of damages that could be caused by Persians to Greeks and vice versa. The third scenario (80% damage by Greeks to Persians and 20% damage by Persians to Greeks) was the most realistic version, confirmed by historical texts.

## 10. References

[1] Osborne, M. J. & Rubinstein, A. (1994) A Course in Game Theory, Cambridge, MA: MIT, ISBN: 9780262150415.

[2] Archan M., Sagar Ch. (2020) Deciphering chaos in evolutionary games. Chaos 30, 121104. doi: 10.1063/5.0029480

[3] Toupo, Danielle F. P., Strogatz, Steven H., Cohen, Jonathan D., Rand, David G. (2015) Evolutionary game dynamics of controlled and automatic decision-making. Chaos: An Interdisciplinary Journal of Nonlinear Science, 25(7), doi:10.1063/1.4927488

[4] Altman E., Gaillard J., Haddad M., Wiecek P. (2012) Dynamic Hawk and Dove Games within Flocks of Birds. In: Hart E., Timmis J., Mitchell P., Nakamo T., Dabiri F. (eds) Bio-Inspired Models of Networks, Information, and Computing Systems. BIONETICS 2011. Lecture Notes of the Institute for Computer Sciences, Social Informatics and Telecommunications Engineering, vol 103. Springer, Berlin, Heidelberg. https://doi.org/10.1007/978-3-642-32711-7_9

[5] Eduardo Souza de Cursi, Uncertainty quantification in game theory, Chaos, Solitons & Fractals, Volume 143, 2021, 110558, ISSN 0960-0779, https://doi.org/10.1016/j.chaos.2020.110558

[6] Benndorf, V, Martínez-Martínez, I, Normann H, 2016, Equilibrium selection with coupled populations in hawk–dove games: Theory and experiment in continuous time, Journal of Economic Theory, Vol. 165, p. 472-486

[7] Smith, M., J., & Price, G. R. (1973) The logic of animal conflict. Nature 246(5427): 15–18.

[8] Lykke, A. F. Jr. (1989) Defining Military Strategy, Military Review 69, no. 5.

[9] Lorenz, E. (1964) The problem of deducing the climate from the governing equations, Tellus 16.

[10] Malthus, T. (1798) An Essay on the Principle of Population, Publisher: J. Johnson, London.

[11] Goodwin, W. (1906) The battle of Salamis. Harvard Studies in Classical Philology, Vol. 17, pp. 74-101.

[12] Green. P. (1998) The Greco – Persian Wars, University of California Press.

[13] Burn, A. R. (1962) Persia and the Greeks, New York: Minerva Press.

[14] Shepherd, W. (2010) Salamis 480 BC: The naval campaign that saved Greece, Osprey Publishing, ISBN: 9781846036842.

[15] Zerefos, C., Solomos, S., Melas, D., Kapsomenakis, J., Repapis, C. (2020) The role of weather during the Greek-Persian "Naval Battle of Salamis" in 480 B.C., Atmosphere, 11(8), 838.

[16] Kyriazis, N., Zouboulakis, M. (2004) Democracy, Sea Power and Institutional Change: An Economic Analysis of the Athenian Naval Law. European Journal of Law and Economics, 17: 117–132, Kluwer Academic Publishers.

[17] Krasanakis, A. (2008) Naval History of the Greek Nation. Athens Press, ISBN: 9789604841592.

### Software

Diks, C., Hommes, C., Panchenko, V., van der Weide, R., "E&F Chaos: A user-friendly software package for nonlinear economic dynamics", Computational Economics, vol. 32, pp. 221- 244, (2008).

McKelvey, Richard D., McLennan, Andrew M., and Turocy, Theodore L., Gambit: Software Tools for Game Theory, Version 16.0.1., (2014). http://www.gambit-project.org.

https://sourceforge.net/projects/maxima/files/Maxima-Windows/5.39.0-Windows/